\documentclass{amsart}
\usepackage{amssymb, amsmath, amsthm, graphics, comment, xspace, enumerate}
\theoremstyle{definition}

\theoremstyle{remark}

\numberwithin{equation}{section}

\begin{document}
\title[ON FREQUENTLY HYPERCYCLIC ABSTRACT HIGHER-ORDER DIFFERENTIAL EQUATIONS]{ON FREQUENTLY HYPERCYCLIC ABSTRACT HIGHER-ORDER DIFFERENTIAL EQUATIONS}

\author{Belkacem Chaouchi}
\address{Lab. de l'Energie et des Syst\' emes Intelligents, Khemis Miliana University, 44225 Khemis Miliana, Algeria}
\email{chaouchicukm@gmail.com}

\author{Marko Kosti\' c}
\address{Faculty of Technical Sciences,
University of Novi Sad,
Trg D. Obradovi\' ca 6, 21125 Novi Sad, Serbia}
\email{marco.s@verat.net}

{\renewcommand{\thefootnote}{} \footnote{
\\ \text{  }  \ \    {\it Key words and phrases.} frequently hypercyclic operators, $C$-regularized semigroups, abstract higher-order differential equations.
\\  \text{  }  \ \ The second named author is partially supported by grant 174024 of Ministry
of Science and Technological Development, Republic of Serbia.}}

\begin{abstract}
In this note, we analyze frequently hypercyclic solutions of abstract higher-order differential equations in separable infinite-dimensional complex Banach spaces.
We essentially apply results from the theory of $C$-regularized semigroups, providing 
several illustrative examples and possible applications.
\end{abstract}

\maketitle

\section{Introduction and Preliminaries}\label{sec:section1}

As it is well-known, the class of frequently hypercyclic linear continuous operators on separable Fr\'echet spaces was introduced by F. Bayart and S. Grivaux in 2006 (\cite{bay1}). Frequent hypercyclicity and various generalizations of this concept are very active fields of research of
a great number of mathematicians working in the field of linear topological dynamics (for more details, we may refer e.g. to
\cite{Baya1}-\cite{boni}, \cite{Grosse} and references cited therein).

Frequently hypercyclic properties of abstract first order differential equations have been studied by E. M. Mangino, A. Peris \cite{man-peris} and E. M. Mangino, M. Murillo-Arcila \cite{man-marina}, within the framework of theory of strongly continuous semigroups, and the second named author \cite{tajwan}, within the theory of integrated and $C$-regularized semigroups. Frequently hypercyclic abstract second order differential equations have been recently investigated in \cite{prcko-fuck it} by using the general notion of $C$-distribution cosine functions and integrated $C$-cosine functions.
Up to now, we do not have any relevant reference treating the operator theoretical aspects of frequently hypercyclic abstract higher-order differential equations. This fact has strongly influenced us to write this paper.

The organization, main ideas and novelties of paper are briefly described as follows. Let $(E,\|\cdot \|)$ be a separable infinite-dimensional complex Banach space. We analyze frequently hypercyclic properties of solutions of the
abstract Cauchy problem\index{abstract Cauchy problem $(ACP_{n})$}
\[(ACP_{n}):\left\{
\begin{array}{l}
u^{(n)}(t)+ A_{n-1}u^{(n-1)}(t)+\cdot \cdot \cdot + A_{1}u^{\prime}(t)+A_{0}u(t)=0,
\ t \geq 0, \\[0.1cm]
u^{(k)}(0)=u_k,\  k=0,\cdot \cdot \cdot, n -1,
\end{array}
\right.
\]
where $A_{0},\cdot \cdot \cdot, A_{n-1}
$ are closed linear operators on $E$ and $u_{0},\cdot \cdot \cdot , u_{n-1}\in E;$ by a strong solution of $(ACP_{n}),$ we mean any $n$-times continuously differentiable function $t\mapsto u(t),$ $t\geq 0$ such that the mappings $t\mapsto A_{i}u^{(i)}(t),$ $t\geq 0$ are continuous for $0\leq i\leq n-1$ and the initial conditions are satisfied (for more details about the well-posedness of $(ACP_{n})$, the reader may consult the monographs 
\cite{x263} by T.-J. Xiao, J. Liang and \cite{knjigaho} by the author).
In order to investigate frequently hypercyclic properties of solutions to $(ACP_{n}),$ we convert this problem into corresponding abstract first order differential equation with appropriately chosen operator matrix acting on product space $E^{n}.$ The proofs of our structural results lean heavily on the use of Lemma 1.1, proved recently in \cite{tajwan}, where we have  considered frequent hypercyclicity for $C$-regularized semigroups following the approach of S. El Mourchid \cite[Theorem 2.1]{samir} and 
E. M. Mangino, A. Peris \cite[Corollary 2.3]{man-peris}. In contrast to the recent research studies of J. A. Conejero, C. Lizama et al. \cite{pako}-\cite{pako2}, where the authors have studied the hypercyclic and chaotic solutions of certain kinds of abstract second and third order differential equations in the spaces of Herzog analytic functions by employing, primarily, the Desch-Schappacher-Webb criterion \cite{fund}, the operator matrix under our consideration is not bounded and as such does not generate a strongly continuous semigrop on $E^{n}$ a priori. This is the main reason why we use the theory of $C$-regularized semigroups in this paper. We construct solutions of $(ACP_{n})$ for initial values $(u_{0},\cdot \cdot \cdot,u_{n-1})$ belonging to a certain proper subspace $\tilde{E}\subseteq E^{n}$ and after that analyze their frequently hypercyclic properties by applying essentially Lemma 1.1, as mentioned above. Motivated by our recent researches \cite{azjm} and \cite{prcko-fuck it}, in Definition 1.1 we introduce the notion of a $(W,\tilde{E},{\mathcal E})$-frequent hypercyclicity. The main goal of Theorem 2.1 is to analyze $(W,\tilde{E},{\mathcal E})$-frequently hypercyclic solutions of some special classes of problems
$(ACP_{n})$ in the case that the operator matrix $p(A)$ obtained after the usual convertion generates an entire $C$-regularized group. After that, we revisit once more the fundamental result \cite[Theorem 5]{neubrander} of F. Neubrander. 
We introduce the notion of $(W,\tilde{E},{\mathcal E}, (D(A_{n-1}))^{n})$-frequent hypercyclicity (Definition 2.1) and   consider $(W,\tilde{E},{\mathcal E}, (D(A_{n-1}))^{n})$--frequently hypercyclic
solutions of $(ACP_{n})$ (Theorem 2.2), provided that the operator $-A_{n-1}$ is the generator of a strongly continuous semigroup on $E$ as well as $D(A_{n-1})\subseteq D(A_{j})$ for $0\leq j\leq n-2$ (cf. also \cite[Theorem 2.10.45]{knjigaho} for a generalization of the above-mentioned theorem to abstract time-fractional differential equations). In our approach, we almost always face the situation $\tilde{E} \neq E^{n},$ which indicates a certain type of subspace frequent hypercyclicity of constructed solutions to $(ACP_{n})$ (in \cite{pako}-\cite{pako2}, the situation in which $\tilde{E} = E^{n}$ can really occur). At the end of paper, we provide several examples and applications of our results. 

Before explaining the notation used, we would like to note that we will not discuss here frequently hypercyclic properties of systems of evolution equations by using the theory of operator matrices developed by K.-J. Engel and his collaborators (for more details about this subject, we refer the reader to the monograph \cite{engel-mb}). 
It would be very tempting to analyze frequently hypercyclic properties of abstract higher-order differential equations by using some other operator theoretical approach which do not use some kind of convertion of the initial equation into a first order matricial equation on product space.

By $L(E)$ we denote the space consisting of all continuous linear mappings from $E$ into $E.$ 
We always assume henceforth that $C\in L(E)$ and $C$ is injective. 
Let $A$ be a closed linear operator with domain $D(A)$ and range $R(A)$ contained in $E,$ and let $CA \subseteq AC.$ Set $D_{\infty}(A):=\bigcap_{k\in {\mathbb N}}D(A^{k}).$ The part of $A$ in a linear subspace ${\tilde E}$ of $E$, $A_{|{\tilde E}}$ shortly, is defined through $A_{|{\tilde E}}:=\{(x,y) \in A : x,\ y\in {\tilde E}\}$ (we will identify an operator and its graph henceforth).
Recall that the $C$-resolvent set of $A,$
denoted by $\rho_{C}(A),$ is defined by
$$
\rho_{C}(A):=\Bigl\{\lambda \in {\mathbb C} : \lambda -A \mbox{ is
injective and } (\lambda-A)^{-1}C\in L(E)\Bigr\}.
$$
In our framework, the
$C$-resolvent set of $A$ consists of those complex numbers $\lambda$
for which the operator $\lambda -A$ is injective and $R(C)
\subseteq R(\lambda -A).$ The resolvent set of $A,$ denoted by $\rho(A),$ is obtained by plugging $C=I.$ 
For every $\lambda \in \rho(A)$ and $n\in {\mathbb N},$ we have that $(D(A^{n}),\| \cdot \|_{n})$ is 
a Banach space, where $\|x\|_{n}:=\sum_{i=0}^{n}\|A^{i}x\|$ ($x\in D(A^{n})$). We denote this space simply by $[D(A^{n})].$
All operator families considered in this paper will be non-degenerate. Set ${\mathbb N}_{n}:=\{1,\cdot \cdot \cdot,n\}$ and ${\mathbb N}_{n}^{0}:={\mathbb N}_{n} \cup \{0\}$ ($n\in {\mathbb N}$).

Suppose that $T\subseteq {\mathbb N}.$
The lower density of $T,$ denoted by $\underline{d}(T),$ is defined through:
$$
\underline{d}(T):=\liminf_{n\rightarrow \infty}\frac{|T \cap [1,n]|}{n}.
$$
If $T\subseteq [0,\infty)$, then the lower density of $T,$ denoted by $\underline{d}(T),$ is defined through:
$$
\underline{d}_{c}(T):=\liminf_{t\rightarrow \infty}\frac{m(T \cap [0,t])}{t},
$$
where $m(\cdot)$ denotes the Lebesgue measure on $[0,\infty).$ 
A linear operator $A$ on $E$ is said to be frequently hypercyclic iff there exists an element $x\in D_{\infty}(A)$ (frequently hypercyclic vector of $A$) such that for each open non-empty subset $V$ of $E$ the set
$\{ n\in {\mathbb N} : A^{n}x\in V \}$ has positive lower density.

Motivated by our recent research study of ${\mathcal D}$-hypercyclic and ${\mathcal D}$-topologically mixing properties of abstract degenerate
Cauchy problems with Caputo fractional derivatives \cite{azjm}, we introduce the following definition (since we are primarily concerned with applications of $C$-regularized semigroups, we will consider only non-degenerate differential equations henceforth; the analysis of frequently hypercyclic abstract time-fractional differential equations is far from being trivial and nothing has been said about this theme so far):

\noindent {\bf Definition 1.1.}  (cf. also \cite[Definition 2]{azjm}) 
Suppose that $\emptyset \neq W\subseteq {\mathbb N}_{n-1}^{0},$ $\tilde{E}$ is a linear subspace of $E^{n}$ and ${\mathcal E}:=(E_{i} : i\in W)$ is a tuple of linear subspaces of $E.$ Then we say that the abstract Cauchy problem $(ACP_{n})$ is $(W,\tilde{E},{\mathcal E})$-frequently hypercyclic iff there exists a strong solution $t\mapsto u(t),$ $t\geq 0$ of $(ACP_{n})$ with the initial values $(u_{0},\cdot \cdot \cdot , u_{n-1})\in \tilde{E}$ satisfying additionally that, for every tuple of open non-empty subsets 
${\mathrm V}:=(V_{i} : i\in W)$ of $E,$ the set
$\bigcap_{i\in W}\{t\geq 0 : u^{(i)}(t)\in V_{i} \cap E_{i}\}$ has positive lower density.
\\

Introduction of Definition 1.1 is also motivated by some recent results about frequently hypercyclic properties of abstract second order differential equations (\cite{prcko-fuck it}). Speaking-matter-of-factly, if the assumptions of \cite[Theorem 1]{prcko-fuck it} are satisfied, then there exists a closed linear subspace $\tilde{E}$ of $E^{2}$ such that the abstract Cauchy problem $(ACP_{2})$ with $A_{1}\equiv 0$ and $A_{0}\equiv -A$ is $(\{0,1\},\tilde{E},{\mathcal E})$-frequently hypercyclic, where ${\mathcal E}=(\pi_{1}(\tilde{E}),\pi_{2}(\tilde{E}))$ and $\pi_{1}(\cdot),\ \pi_{2}(\cdot)$ denote the first and second projection, respectively. It is also worth noting that the spectral conditions of \cite[Theorem 1]{prcko-fuck it} are particularly satisfied for a substantially large class of abstract incomplete second order differential equations.

We will use the following definition:

\noindent {\bf Definition 1.2.} 
Let $A$ be a closed linear operator. If
there exists a strongly continuous operator family
$(T(t))_{t\geq 0}\subseteq L(E)$ such that:
\begin{itemize}
\item[(i)] $T(t)A\subseteq AT(t)$, $t\geq 0$,
\item[(ii)] $T(t)C=CT(t)$, $t\geq 0$,
\item[(iii)] for all $x\in E$ and $t\geq 0$: $\int_{0}^{t}T(s)x\,ds\in D(A)$ and
\[
A\int\limits_{0}^{t}T(s)x\,ds=T(t)x-Cx,
\]
\end{itemize}
then it is said that $A$ is a subgenerator of a (global)
$C$-regularized semigroup $(T(t))_{t\geq 0}$.
\\

It is well-known that $T(t)T(s)=T(t+s)C$ for all $t,\ s\geq 0.$ The integral generator of $(T(t))_{t\geq 0}$ is defined by
\[
\hat{A}:=\Biggl\{(x,y)\in E\times
E:T(t)x-Cx=\int\limits^{t}_{0}T(s)y\,ds,\;t\geq
0\Biggr\}.
\]
We know that the integral generator of $(T(t))_{t\geq 0}$
is a closed linear operator which is an
extension of any subgenerator of $(T(t))_{t\geq 0}$
and satisfies $\hat{A}=C^{-1}AC$ for any subgenerator $A$ of $(T(t))_{t\geq 0}.$ If for each fixed element $x\in E$ the mapping $t\mapsto T(t)x,$ $t\geq 0$ can be extended to an entire function, then we say that $(T(t))_{t\geq 0}$ is an entire $C$-regularized group with subgenerator $A$ and integral generator $\hat{A}$ (\cite{l1}). Furthermore, it is said that $(T(t))_{t\geq 0}$ is frequently hypercyclic iff there exists an element $x\in E$ (frequently hypercyclic vector of $(T(t))_{t\geq 0}$) such that the mapping $t\mapsto C^{-1}T(t)x,$ $t\geq 0$ is well-defined, continuous and that for each open non-empty subset $V$ of $E$ the set $\{ t\geq 0 : C^{-1}T(t)x\in V \}$ has positive lower density (\cite{tajwan}).

Throughout the whole paper, we will essentially employ the following result (see \cite[Theorem 2.10(ii)]{tajwan}):

\noindent {\bf Lemma 1.1.}  {\it 
Let $t_0>0$ and let $A$ be a subgenerator of a global $C$-regularized semigroup $(S_{0}(t))_{t\geq 0}$ on $E.$ Suppose that $R(C)$ is dense in $E.$ Set $T(t)x:=C^{-1}S_{0}(t)x,$ $t\geq 0,$ $x\in Z_{1}(A).$
Suppose, further, that there exists a family $(f_{j})_{j\in \Gamma}$ of twice continuously differentiable mappings $f_{j} : I_{j} \rightarrow E$ such that $I_{j}$ is an interval in ${\mathbb R}$ and $Af_{j}(t) = itf_{j}(t)$ for every
$t \in I_{j} ,$ $ j \in \Gamma .$ 
Set $\tilde{E}:=\overline{span\{f_{j}(t) : j \in \Gamma,\ t \in I_{j}\}}$.
Then $A_{|\tilde{E}}$ is a subgenerator of a global $C_{|\tilde{E}}$-regularized semigroup $(S_{0}(t)_{| \tilde{E}})_{t\geq 0}$ on $\tilde{E},$
$(S_{0}(t)_{|\tilde{E}})_{t\geq 0}$ is frequently hypercyclic in $\tilde{E}$ and the operator $T(t_{0})_{|\tilde{E}}$
is frequently hypercyclic in $\tilde{E}.$} \\

For more details about $C$-regularized 
semigroups and their applications, we refer the reader to the monographs \cite{l1} by R. deLaubenfels and \cite{knjigah}-\cite{knjigaho} by the author.  

\section{Formulation and Proof of Main Results}\label{idiot}
 
In the formulation of our first structural result, we assume that 
$N,\ n\in {\mathbb N}$ and $iA_{j},\ 1\leq j\leq N$ are commuting
generators of bounded $C_{0}$-groups on $E$ (here, $i$ denotes the imaginary unit). Define $A:=(A_{1},\cdot
\cdot \cdot, A_{N})$ and $A^{\eta}:=A_{1}^{\eta_{1}}\cdot \cdot
\cdot A_{N}^{\eta_{N}}$ for any $\eta=(\eta_{1},\cdot \cdot \cdot,
\eta_{N})\in {{\mathbb N}_{0}^{N}}.$ If $P(\xi)=[p_{ij}(\xi)]_{n\times n}$ is an arbitrary matrix of complex polynomials in variable $\xi \in {\mathbb R}
^{N},$ then we can write $P(\xi)=\sum_{|\eta |\leq m}P_{\eta}\xi^{\eta}$ for a certain integer $m\in {\mathbb N}$ and for certain complex matrices $P_{\eta}$ of format $n\times n.$ 
We know that the operator $P(A) :=\sum_{|\eta |\leq m}P_{\eta}A^{\eta}$ acting with its maximal domain is closable on $E^{n};$ moreover, the following holds:

\noindent {\bf Lemma 2.1.} (\cite{l1}, \cite{knjigaho}) {\it
There exists an injective operator $C\in L(E^{n})$ with dense range in $E^{n}$ such that 
the operator $\overline{P(A)}$ generates an entire $C$-regularized group $(T(t))_{t\geq 0}$ on $E^{n}$ such that $T(t)\vec{x}\in D_{\infty}(P(A))$ for all $\vec{x}\in E^{n}.$} \\

Let $\pi_{j} : E^{n}\rightarrow E$ be the $j$-th projection ($1\leq j\leq n$), let $p_{0}(\xi),\cdot \cdot \cdot,p_{n-1}(\xi)$ be complex polynomials in variable $\xi \in {\mathbb R}^{N},$ and let   
$$
p(A):=\left[ \begin{array}{ccccc}
0 & I & 0  & \cdot \cdot \cdot & 0 \\
 0 & 0 &  I & \cdot \cdot \cdot & 0  \\
\cdot & \cdot &  \cdot & \cdot \cdot \cdot & \cdot \\
0 & 0 & 0 &  \cdot \cdot \cdot  & I \\
-A_{0} & -A_{1} & -A_{2} & \cdot \cdot \cdot & -A_{n-1}  \end{array} \right],
$$
where $A_{i}:=p_{i}(A)$ for $0\leq i\leq n-1.$
Then we have the following:

\noindent {\bf Theorem 2.1.} 
{\it Suppose that there exists a family $(F_{j})_{j\in \Gamma}$ of twice continuously differentiable mappings $F_{j} : I_{j} \rightarrow E^{n}$ such that $I_{j}$ is an interval in ${\mathbb R}$ and $p(A)F_{j}(t) = itF_{j}(t)$ for every
$t \in I_{j} ,$ $ j \in \Gamma .$ 
Set $\tilde{E}:=\overline{span\{F_{j}(t) : j \in \Gamma,\ t \in I_{j}\}}$. Then the abstract Cauchy problem $(ACP_{n})$ is $({\mathbb N}_{n-1}^{0},\tilde{E}, {\mathcal E})$-frequently hypercyclic with
${\mathcal E}:=(\pi_{1}(\tilde{E}),\cdot \cdot \cdot, \pi_{n}(\tilde{E})).$}

\begin{proof}
By Lemma 2.1, we know that there exists an injective operator $C\in L(E^{n})$ with dense range in $E^{n}$ such that 
the operator $\overline{P(A)}$ generates an entire $C$-regularized group $(T(t))_{t\geq 0}$ on $E^{n}$ such that $T(t)\vec{x}\in D_{\infty}(P(A))$ for all $\vec{x}\in E^{n}.$
Furthermore, the injective operator $C$ can be chosen such that $C(span\{F_{j}(t) : j \in \Gamma,\ t \in I_{j}\})=span\{F_{j}(t) : j \in \Gamma,\ t \in I_{j}\}$
and that $C(\tilde{E})$ is a dense linear subspace of $\tilde{E};$ see \cite[Remark 14(ii)]{cds}.
Due to Lemma 1.1, we know that $(S_{0}(t)_{|\tilde{E}})_{t\geq 0}$ is frequently hypercyclic in $\tilde{E},$ which implies that there exists a vector $\vec{x}\in \tilde{E}$ such that for each open non-empty subset $V$ in $E^{n}$ the set $\{t\geq 0 : C^{-1}T(t)\vec{x} \in \tilde{E} \cap V\}$ has positive lower density. 
Since $C(\tilde{E})$ is a dense linear subspace of $\tilde{E},$ it readily follows that for each open non-empty subset $V$ in $E^{n}$ the set $\{t\geq 0 : T(t)\vec{x} \in \tilde{E} \cap V\}$ has positive lower density, as well. On the other hand, the function $t\mapsto T(t)\vec{x},$ $t\geq 0$ is a unique solution of the abstract Cauchy problem
$\vec{U}^{\prime}(t)=\overline{p(A)}\vec{U}(t),$ $t\geq 0;$ $\vec{U}(0)=C\vec{x}.$
Furthermore, $T(t)\vec{x}\in D_{\infty}(P(A))$ so that the function $t\mapsto T(t)\vec{x},$ $t\geq 0$ is a unique solution of the abstract Cauchy problem
$\vec{U}^{\prime}(t)=p(A)\vec{U}(t),$ $t\geq 0;$ $\vec{U}(0)=C\vec{x},$ actually. It is clear that the first, second,..., the $n$-th component of $T(\cdot)\vec{x}$ is a unique solution of $(ACP_{n}),$ its first derivative,..., its $(n-1)$-derivative, respectively, with the initial conditions $u_{j}=\pi_{j+1}(C\vec{x}),$ $0\leq j\leq n-1.$ This simply implies the required conclusion.    
\end{proof}

\noindent {\bf Remark 2.1.} 
The most important case for applications is $N=1.$
In this case, let us assume that $f_{j} : I_{j} \rightarrow E$ is a twice continuously differentiable mapping,
$g_{j} : \{it \,  ; \,  t\in I_{j}\}\rightarrow {\mathbb C} \setminus \{0\}$ is a scalar-valued mapping and
$Af_{j}(t)=g_{j}(it)f_{j}(t),$ $t \in I_{j}$ ($j\in \Gamma$). If
\begin{equation}\label{prcko}
(it)^{n}+\sum_{l=0}^{n-1}(it)^{l}P_{l}\bigl(g_{j}(it)\bigr)=0,\quad t\in I_{j},\ j\in \Gamma ,
\end{equation}
then the assumptions of Theorem 2.1 are satisfied with 
\[
F_{j}(t):=\bigl[f_{j}(t)\ itf_{j}(t) \cdot \cdot \cdot \ (it)^{n-1}f_{j}(t) \bigr]^{T},\quad t\in I_{j},\ j\in \Gamma;
\]
see e.g. \cite[Example 1(ii)]{azjm}. 
\\

We continue by observing that Definition 1.1 does not enable one to thoroughly investigate frequently hypercyclic solutions of some important classes of abstract higher-order differential equations already examined in the existing literature. For example, F. Neubrander has analyzed in \cite{neubrander} the well-posedness results for $(ACP_{n})$ by reduction this problem into a first order matricial system, employing the matrix
$$
\Delta:=\left[ \begin{array}{ccccc}
-A_{n-1} & I & 0  & \cdot \cdot \cdot & 0 \\
-A_{n-2} & 0 &  I & \cdot \cdot \cdot & 0  \\
\cdot & \cdot &  \cdot & \cdot \cdot \cdot & \cdot \\
-A_{1} & 0 & 0 &  \cdot \cdot \cdot  & I \\
-A_{0} & 0 & 0 & \cdot \cdot \cdot & 0  \end{array} \right].
$$
The operator matrix
$$
\Psi:=\left[ \begin{array}{cccccc}
I & 0 & 0  & \cdot \cdot \cdot & 0 & 0 \\
-A_{n-1} & I & 0 & \cdot \cdot \cdot & 0 & 0 \\
-A_{n-2} & -A_{n-1} &  I & \cdot & \cdot \cdot \cdot &  \\
\cdot & \cdot & \cdot &  \cdot \cdot \cdot  & I & 0 \\
-A_{1} & -A_{2} & -A_{3} & \cdot \cdot \cdot & -A_{n-1}  & I\end{array} \right]
$$
plays an important role in his analysis, as well. 

We will use the following notion:

\noindent {\bf Definition 2.1.} 
Suppose that $\lambda \in \rho(\Delta),$ $\emptyset \neq W\subseteq {\mathbb N}_{n-1}^{0},$ $\tilde{E}$ is a linear subspace of $(D(A_{n-1}))^{n}$ and ${\mathcal E}:=(E_{i} : i\in W)$ is a tuple of linear subspaces of $E.$ Then we say that the abstract Cauchy problem $(ACP_{n})$ is $(W,\tilde{E},{\mathcal E}, (D(A_{n-1}))^{n})$-frequently hypercyclic iff there exists a strong solution $t\mapsto u(t),$ $t\geq 0$ of $(ACP_{n})$ with the initial values $(u_{0},\cdot \cdot \cdot , u_{n-1})\in \tilde{E}$ satisfying additionally that, for every open non-empty subset
${\mathrm V}$ of $E^{n},$ the set
$\bigcap_{i\in W}\{t\geq 0 : u^{(i)}(t)+\sum_{j=1}^{i}A_{n-i}u^{(i-j)}(t) \in \pi_{i+1}((\lambda-\Delta)^{-n}({\mathrm V})) \cap E_{i}\}$ has positive lower density.\\

This definition is a good one and does not depend on the choice of number $\lambda \in \rho(\Delta).$ This follows from the fact that for each $\lambda \in \rho(\Delta)$ the mapping $\Pi : E^{n} \rightarrow [D(\Delta^{n})]$ given by $\Pi \vec{x}:=(\lambda -\Delta)^{-n}\vec{x},$ $\vec{x}\in E^{n}$ is a linear topological isomorphism so that $\{(\lambda-\Delta)^{-n}({\mathrm V}): {\mathrm V}\mbox{ is an open non-empty subset of }E^{n}\}$ is equal to the set of all open non-empty subsets of $[D(\Delta^{n})]$ and therefore independent of $\lambda \in \rho(\Delta).$

Our second structural result reads as follows:

\noindent {\bf Theorem 2.2.} 
{\it Suppose that the operator $-A_{n-1}$ is the generator of a strongly continuous semigroup on $E$ as well as $D(A_{n-1})\subseteq D(A_{j})$ for $0\leq j\leq n-2.$ Suppose, further, that there exists a family $(F_{j})_{j\in \Gamma}$ of twice continuously differentiable mappings $F_{j} : I_{j} \rightarrow E^{n}$ such that $I_{j}$ is an interval in ${\mathbb R}$ and $\Delta F_{j}(t) = itF_{j}(t)$ for every
$t \in I_{j} ,$ $ j \in \Gamma .$ 
Set 
\begin{equation}\label{prevaz}
\tilde{E}:=\overline{span\bigl\{F_{j}(t) : j \in \Gamma,\ t \in I_{j}\bigr\}}^{[D(\Delta^{n})]}.
\end{equation} 
Then the abstract Cauchy problem $(ACP_{n})$ is $({\mathbb N}_{n-1}^{0},\Psi^{-1}(\tilde{E}), {\mathcal E})$-frequently hypercyclic with
${\mathcal E}:=(\pi_{1}(\tilde{E}),\cdot \cdot \cdot, \pi_{n}(\tilde{E})).$}\\

\begin{proof}
By the proof of \cite[Theorem 5]{neubrander}, we know the following:
\begin{itemize}
\item[(i)] The operator $\Delta$ generates a strongly continuous semigroup $(T(t))_{t\geq 0}$ on $E^{n}$ and therefore there exists $\lambda \in \rho(\Delta).$ 
\item[(ii)] The mapping $\Psi$ is a bijection between the spaces $(D(A_{n-1}))^{n}$ and $D(\Delta^{n}).$
\item[(iii)] For every $\vec{x}\in D(\Delta^{n}),$ the mapping $t\mapsto \pi_{1}(T(t)\vec{x}),$ $t\geq 0$ is a strong solution of problem 
$(ACP_{n})$ with the initial value $\vec{y}=\Psi^{-1}\vec{x}\in (D(A_{n-1}))^{n}.$
\end{itemize}
From (i), we may deduce that the operator $\Delta_{|D(\Delta^{n})}$ generates a strongly continuous semigroup on the space $[D(\Delta^{n})];$ see e.g. \cite[Chapter II.5]{engel}. By
Lemma 1.1, it follows that there exists a vector $\vec{x}\in D(\Delta^{n}) \cap \tilde{E}$ such that 
for every open non-empty subset ${\mathrm V}'$ in $[D(\Delta^{n})],$
the set $\{t\geq 0 : T(t)\vec{x} \in {\mathrm V}' \cap \tilde{E}\}$ has positive lower density. Since $\{t\geq 0 : T(t)\vec{x} \in {\mathrm V}' \cap \tilde{E}\}\subseteq \bigcap_{i=1}^{n}\{ t\geq 0 : \pi_{i}(T(t)\vec{x})\in \pi_{i}( {\mathrm V}') \cap \pi_{i}(\tilde{E})\},$ the required assertion follows from a simple analysis involving (i)-(iii) and the fact that the $(i+1)$-projection of $T(\cdot)\vec{x}$ equals $ u^{(i)}(\cdot)+\sum_{j=1}^{i}A_{n-i}u^{(i-j)}(\cdot)$ for $0\leq i\leq n-1,$ where $u(\cdot):=\pi_{1}(T(\cdot)\vec{x})$ is a unique strong solution of problem $(ACP_{n})$ with initial value $y=\Psi^{-1}\vec{x}$ (see the equation \cite[(1), p. 267]{neubrander}).
\end{proof}

\noindent {\bf Remark 2.2.} 
Let us assume that 
$P_{0},\cdot \cdot \cdot, P_{n-1}$ are complex polynomials in one variable,
$f_{j} : I_{j} \rightarrow E$ is a twice continuously differentiable mapping, 
$g_{j} : \{it \,  ; \,  t\in I_{j}\}\rightarrow {\mathbb C} \setminus \{0\}$ is a scalar-valued mapping and
$Af_{j}(t)=g_{j}(it)f_{j}(t),$ $t \in I_{j}$ ($j\in \Gamma $). If (\ref{prcko}) holds,
then the assumptions of Theorem 2.2 are satisfied with $A_{s}:=P_{s}(A)$ ($0\leq s \leq n-1$) and
\begin{equation}\label{ajsa}
F_{j}(t):=\bigl[F_{j1}(t)\ F_{j2}(t) \cdot \cdot \cdot \ F_{jn}(t) \bigr]^{T},\quad t\in I_{j},\ j\in \Gamma,
\end{equation}
where, for $1\leq s\leq n,$ 
\begin{equation}\label{prckoled}
F_{js}(t):=\sum\limits_{l=0}^{s-2}(it)^{l}A_{n-s+1+l}f_{j}(t)+(it)^{s-1}f_{j}(t),\quad t\in I_{1},\ j\in \Gamma.
\end{equation}

It is worth noting that Theorem 2.1 and Theorem 2.2 provide also sufficient spectral conditions for certain types of (subspace) topologically mixing properties and (subspace) Devaney chaoticity of solutions to $(ACP_{n});$ see \cite{azjm} for more details.

We close the paper by providing some illustrative examples and applications.

\noindent {\bf Example 2.1.} (\cite{{coni-man}}, \cite{metafune})
Suppose that $E:=L^{2}({\mathbb R}),$ $c>\frac{b}{2}>0,$ the operator ${\mathcal A}_{c}$ is defined by $D({\mathcal A}_{c}):=\{u\in L^{2}({\mathbb R})\ \cap \
W^{2,2}_{loc}({\mathbb R}) : {\mathcal A}_{c}u\in L^{2}({\mathbb
R})\},$ 
${\mathcal A}_{c}u:=u^{\prime
\prime}+bxu^{\prime}+cu,$ $u\in D({\mathcal A}_{c}),$ $\Omega :=\{ \lambda \in {\mathbb C} : \Re \lambda <c-\frac{b}{2}\},$
$f_{1}(\lambda):={\mathcal
F}^{-1}(e^{-\frac{\xi^{2}}{2b}}\xi
|\xi|^{-(2+\frac{\lambda-c}{b})})(\cdot),$ $\lambda \in \Omega $ and $f_{2}(\lambda):={\mathcal F}^{-1}(e^{-\frac{\xi^{2}}{2b}}
|\xi|^{-(1+\frac{\lambda-c}{b})})(\cdot),$ $\lambda \in \Omega$ (here, ${\mathcal
F}^{-1}$ denotes the inverse Fourier transform on the real line). Consider the equation
\begin{equation}\label{rajko}
u^{(n)}(t)+\sum_{j=1}^{n-1}P_{j}({\mathcal  A}_{c})u^{(j)}(t)=0,\quad
t\geq 0,
\end{equation}
where the operator $-P_{n-1}({\mathcal  A}_{c})$ is the generator of a strongly continuous semigroup on $E$ and $P_{n-1}(\cdot)$ is a non-zero complex polynomial whose degree is less or equal than the degree of any other non-zero complex polynomial $P_{j}(\cdot)$ for $0\leq j\leq n-2.$ Keeping in mind that ${\mathcal A}_{c}f_{1,2}(\lambda)=\lambda f_{1,2}(\lambda),$ $\lambda \in \Omega ,$ the validity of equality
$$
(it)^{n}+\sum\limits_{j=0}^{n-1}(it)^{j}P_{j}(it)=0,\quad t\in {\mathbb R},
$$
implies that the equation (\ref{prcko}) holds with $\Gamma:=\{1,2\},$ $I_{1}=I_{2}:={\mathbb R}$ and $g_{1}(it)=g_{2}(it):=it,$ $t\in {\mathbb R}.$ Therefore, Theorem 2.2 is applicable so that the abstract Cauchy problem (\ref{rajko}), equipped with initial conditions $u^{(j)}(0)=u_{j}$ for $0\leq j\leq n-1,$ is $({\mathbb N}_{n-1}^{0},\Psi^{-1}(\tilde{E}), {\mathcal E})$-frequently hypercyclic, where $F_{j}(\cdot)$ is given by (\ref{ajsa})-(\ref{prckoled}) for $j=1,2,$ $\tilde{E}$ is defined by (\ref{prevaz}) and
${\mathcal E}:=(\pi_{1}(\tilde{E}),\pi_{2}(\tilde{E}),\cdot \cdot \cdot,\pi_{n}(\tilde{E})).$ A concrete example where the above assumptions are satisfied can be plainly constructed in the following way: 
Since the operator $A_{1}:={\mathcal A}-c_{1}$ generates a strongly continuous semigroup for any $c_{1}\in {\mathbb C}, $ we can 
choose $P_{n-1}(z):=-z+c_{1},$ $P_{2}(z):=-c_{1}z$ and $P_{j}(z):=0$ for $j>2.$

\noindent {\bf Example 2.2.} (\cite{ji})
Suppose that $0< \gamma \leq 1,$ $a>0,$ $p>2$ and $X$ is a symmetric
space of non-compact type\index{symmetric spaces of non-compact
type} and rank one.  Then the Laplace-Beltrami operator $-\Delta_{X,p}^{\natural}$ generates a strongly continuous semigroup on $X$ and we
know that $int(P_{p})\subseteq
\sigma_{p}(\Delta_{X,p}^{\natural}),$ where $P_{p}$ denotes the parabolic domain defined in \cite{ji}. Suppose that (2.1) holds with
$\Gamma=\{1\},$ the function $g_{1}(it)=it,$ $t\in I_{1}$ ($I_{1}$ is a suitable chosen subinterval of ${\mathbb R}$) and certain complex polynomials $P_{0}(\cdot),\cdot\cdot \cdot,P_{n-1}(\cdot).$ Then Theorem 2.2 is applicable with operators $A_{l}:=P_{l}(A)$ ($l\in {\mathbb N}_{n-1}^{0}$), provided that the operator $-A_{n-1}$ is the generator of a strongly continuous semigroup on $E$ and $P_{n-1}(\cdot)$ is a non-zero complex polynomial whose degree is less or equal than the degree of any other non-zero complex polynomial $P_{j}(\cdot)$ for $0\leq j\leq n-2.$ 

\noindent {\bf Example 2.3.}  
Let us recall that a measurable function
$\rho : {\mathbb R} \rightarrow (0,\infty)$ is called an
admissible weight function\index{admissible weight function} iff there exist constants $M\geq 1$ and
$\omega \in {\mathbb R}$ such that $\rho(t)\leq Me^{\omega
|t'|}\rho(t+t')$ for all $t,\ t'\in {\mathbb R}.$ For such a function
$\rho(\cdot),$ we consider the following Banach spaces:
$$
L^{p}_{\rho}({\mathbb R}):=\bigl\{ u : {\mathbb R} \rightarrow
{\mathbb C} \ ; u(\cdot) \mbox{ is measurable and } ||u||_{p} <\infty \bigr\},
$$
where $p\in [1,\infty)$ and $||u||_{p}:=(\int _{\mathbb R}|u(t)|^{p}
\rho(t)\, dt)^{1/p},$ as well as
$$
C_{0,\rho}({\mathbb R}):=\bigl\{u : {\mathbb R} \rightarrow {\mathbb
C} \ ; u(\cdot) \mbox{ is continuous and }\lim _{t \rightarrow
\infty}u(t)\rho(t)=0\bigr\},
$$
with $||u||:=\sup _{t\in {\mathbb R}}|u(t)\rho(t)|.$ 
It is well-known that the operator $A:=d/dt$ equipped with domain $D(A):=\{u\in E : u^{\prime}\in E,\ u(\cdot)\mbox{ is absolutely continuous}\}$ generates a strongly continuous translation group on $E$ (see \cite[Lemma 4.6]{fund}). If we assume that, for every $\lambda \in i{\mathbb R},$ the function $t\mapsto e^{\lambda t},$ $t\in {\mathbb R}$ belongs to the space $E$ 
and the equation (2.1) holds with
$\Gamma=\{1\},$ the function $g_{1}(it)=it,$ $t\in I_{1}={\mathbb R}$ and certain complex polynomials $P_{0}(\cdot),\cdot\cdot \cdot,P_{n-1}(\cdot),$ then Theorem 2.1 is applicable with operators $A_{l}:=P_{l}(A)$ ($l\in {\mathbb N}_{n-1}^{0}$) and $N=1,$ where $iA$ is the generator of a bounded $C_{0}$-group on $E.$

 \bigskip

{\small\rm\baselineskip=10pt
 \baselineskip=10pt
 \qquad Belkacem Chaouchi\par
 \qquad Lab. de l'Energie et des Syst\' emes Intelligents\par
 \qquad Khemis Miliana University\par
 \qquad 44225 Khemis Miliana, Algeria\par
 \qquad {\tt chaouchicukm@gmail.com}

 \bigskip \smallskip

 \qquad Marko Kosti\' c\par
 \qquad Faculty of Technical Sciences\par
 \qquad University of Novi Sad\par
 \qquad Trg D. Obradovi\' ca 6 73\par
 \qquad 21125 Novi Sad, Serbia\par
 \qquad {\tt marco.s@verat.net}
 }

 \end{document}